\documentclass[12pt]{article}
\usepackage{amsmath}
\usepackage{amssymb}
\usepackage{geometry}
\usepackage{amscd}
\usepackage{xypic}
\usepackage{mathrsfs}
\usepackage{amsthm}
\usepackage{color}
\usepackage{fancyhdr}
\usepackage{graphicx}
\usepackage{pgf,tikz}
\usetikzlibrary{arrows}
\usetikzlibrary{calc,intersections,through,backgrounds}
\usepackage{array}
\usepackage{float}
\usepackage{hyperref}

\newtheorem{theorem}{Theorem}[section]
\newtheorem{definition}{Definition}[section]
\newtheorem{proposition}{Proposition}[section]

\theoremstyle{definition}

\newtheorem{remark}{Remark}[section]

\begin{document}

\title{Computing Dixmier Invariants and Some Geometric Configurations of Quartic Curves with 2 Involutions 
}

\author{Dun Liang      
}

\maketitle

\begin{abstract} In this paper we consider plane quartics with to involutions. We compute the Dixmier invariants, the bitangents and the Matrix representation problem of these curves, showing that they have symbolic solutions for the last two questions.
\end{abstract}
\section{Introduction}
We consider algebraic varieties over the field $K=\overline{\mathbb Q}$, the algebraic closure of the rational numbers $\mathbb Q$ in the field of complex numbers $\mathbb C$. Let $C$ be a general quartic curve in the projective plane ${\mathbb P}^2_K$ (or ${\mathbb P}^2$ in short) defined by the general equation
\begin{equation}\label{genC}
f(x,y,z)=\sum_{i+j+k=4}a_{ijk}x^iy^jz^k=a_{400}x^4+a_{310}x^3y+\cdots +a_{004}z^4
\end{equation}
where $a_{ijk}$ are the corresponding coefficients for $i,j,k\in \{0,1,2,3,4\}$. According to the theory of algebraic curves, a generic plane quartic is smooth, and the genus $g(C)=3$. The study of such curves is an important topic in clasical algebraic geometry (see \cite{CAG}), and in modern times, people study the moduli space ${\cal M}_3$ of such curves.

In this paper we study smooth plane quartics with ${\mathbb Z}/2\times {\mathbb Z}/2$-automorphisms. Since smooth plane quartics are non-hyperelliptic genus 3 curves, the equation (\ref{genC}) is the canonical model of teh curve $C$ if $C$ is smooth. If $C$ has an automorphism $\phi$, then $\phi$ is a projective linear transformation on ${\mathbb P}^2$ with respect to the equation (\ref{genC}).

The classification of automorphisms of smooth plane quartics is given by \cite{Kantor} and \cite{Wiman}, people can find a full list in many references nowadays, such as \cite{CAG}, \cite{Kohel}. Explicity, smooth quartics with  ${\mathbb Z}/2\times {\mathbb Z}/2$-automorphisms should be isomorphic to one of the following curves:

\begin{table}[h]\label{table}\caption{Curves with  ${\mathbb Z}/2\times {\mathbb Z}/2$-Automorphisms}
\begin{center}\begin{tabular}{|c|c|c|}
\hline 
Name & Automorphism Group & Equation \\ 
\hline 
$X_4$ & ${\mathbb Z}/2\times {\mathbb Z}/2$ & $x^4+y^4+z^4+rx^2y^2+sy^2z^2+uz^2x^2=0$ \\ 
\hline 
$X_{16}$ & $D_8$ & $x^4+y^4+z^4+rx^2y^2+s(y^2z^2+z^2x^2)=0$ \\ 
\hline 
$X_{24}$ & ${\mathfrak S}_4$ & $x^4+y^4+z^4+r(x^2y^2+y^2z^2+z^2x^2)=0$ \\ 
\hline 
$X_{96}$ & $({\mathbb Z}/4\times {\mathbb Z}/4)\ltimes {\mathfrak S}_3$ & $x^4+y^4+z^4=0$ \\ 
\hline 
\end{tabular} \end{center}
\end{table}

We consider three geometric informations about these curves, the Dixmier invariants, the bitangents, and the matrix representation probelm.

Like the $j$-invariant of elliptic curves, the quartics have their own invariants, the Dixmier-Ohno invariants \cite{Dixmier,Elsenhans,Ohno}. In fact Ohno \cite{Ohno} gives covariants. We only compute the Dixmier invariants \cite{Dixmier}. There are $7$ of them
$I_3,I_6,I_9,I_{12},I_{15},I_{18}$ and the discriminant $I_{27}$. We will not compute the discriminant $I_{27}$. In fact all the curves in Table \ref{table} are smooth and the discriminant is the invariant to judge the smoothness. We use Maxima \cite{Maxima} to compute the invariants of the curves in Table \ref{table} in Section 2. The invariants of the general family $X_4$ is symmetric with respect to the parameters, so we write the invariants as polymomials of the elementary symmetric functions (see Proposition \ref{propX4}).

The bitangents of the plane quartics is an important topic in classical algebraic geometry (see \cite{Hesse1,Hesse2,Kantor,Jacobi}). As divisors on the curve, they are related to the theta characteristics of the curve (see \cite{Mumford,Steiner}). We use the idea in \cite{Sturmfels} to compute the bitangents. The main result is Theorem \ref{BitanX4}, showing that all the curves in Table \ref{table} have symbolic solutions of all 28 bitangents. In the proof, we use {\tt Macaulay2} \cite{Macaulay2} to make the elimination, and use Maxima to compute the solutions.

The matrix representation problem (see \cite{Vinnikov1,Helton,Vinnikov2}) asks if a quartic homogeneous polynomial could be written as 
$$\det (xA+yB+zC)$$
for some symmetric matrices $A,B,C$. We use the idea in \cite{Sturmfels2} to compute this problem. We first use Maxima to compute the determinant, and then argue the conditions on the entries. At the end we reach Theorem \ref{MRPX4}, showing that this problem has a symbolic solution for $X_4$.

\section{Dixmier Invariants of $X_4$, $X_{16}$, $X_{24}$ and $X_{96}$}
\subsection{Dixmier Invariants of Plane Quartics}\label{Dix}

Our notations follows from \cite{Kohel}. First we introduce some notations. In general, let $f\in K[x_1,\ldots, x_n]$ be a polynomial, we use $D_f$ to denote the differential operator determined by $f$. Explicitly, let 
\begin{equation}\label{genf} f=f(x_1,\ldots ,x_n)= \sum_{(i_1,\ldots , i_n)\in {\mathbb Z}^n_+} a_{i_1,\ldots ,i_n}x_1^{i_1}\cdots x_n^{i_n}\end{equation}
where $a_{i_1,\ldots ,i_n}\in K$ be the coefficient of the monomial $x_1^{i_1}\cdots x_n^{i_n}$ for $(i_1,\ldots , i_n)\in {\mathbb Z}^n_+$ and (\ref{genf}) be a finite sum. For the rest of this paper, we will not emphasize that the powers $i_1,\ldots , i_n$ are non-negative integers again. 

The map $D_f$ means
$$\begin{array}{cccc}
D_f : & K[x_1,\ldots, x_n] & \longrightarrow & K[x_1,\ldots, x_n] \\ \\ 
& g(x_1,\ldots , x_n) & \longmapsto & \displaystyle{\sum_{(i_1,\ldots , i_n)\in {\mathbb Z}^n_+} a_{i_1,\ldots ,i_n}\frac{\partial^{i_1+\cdots + i_n}}{\partial x_1^{i_1}\cdots{\partial x_n^{i_n}}} g(x_1,\ldots , x_n)}.
\end{array}
$$
If we use $D(f,g)$ to denote $D_f(g)$, $\forall f,g\in K[x_1,\ldots ,x_n]$, then the map 
$$D:K[x_1,\ldots ,x_n] \times K[x_1,\ldots ,x_n] \longrightarrow K[x_1,\ldots ,x_n]$$ has some obvious properties as the following:

\begin{itemize}
\item $D$ is bilinear.
\item Let ${\rm deg}(f)$ be the degree of $f$ for all $f\in K[x_1,\ldots ,x_n]$.
 Let $f,g\in K[x_1,\ldots ,x_n]$. If ${\rm deg}(f)> {\rm deg}(g)$, then $D_f(g)=0$. If 
${\rm deg}(f)> {\rm deg}(g)$, then $D_f(g)\leq {\rm deg}(g)-{\rm deg}(f)$.  Let $f=x_1^{i_1}\cdots x_n^{i_n}$ and $g= x_1^{j_1}\cdots x_n^{j_n}$ be two monomials such that ${\rm deg}(f)= {\rm deg}(g)$, then $D_f(g)=i_1!\cdots i_n!\delta_{fg}$ where $\delta_{fg}$ is the Kronecker delta of $f$ and $g$.
\end{itemize}

For any $f\in K[x_1,\ldots ,x_n]$, let $H(f)$ be the half Hessian matrix of $f$. For example, if $f\in K[x,y,z]$, then 
$$H(f)=\begin{pmatrix}
\displaystyle{\frac{\partial^2}{\partial x^2}} & \displaystyle{\frac{\partial^2}{\partial x \partial y}} & \displaystyle{\frac{\partial^2}{\partial x \partial z}} \\ 
 \displaystyle{\frac{\partial^2}{\partial x \partial y}}& \displaystyle{\frac{\partial^2}{\partial y^2}} & \displaystyle{\frac{\partial^2}{\partial y \partial z}} \\
\displaystyle{\frac{\partial^2}{\partial x \partial z}} & \displaystyle{\frac{\partial^2}{\partial y \partial z}} & \displaystyle{\frac{\partial^2}{\partial z^2}}
\end{pmatrix}.$$
Let $H^*(f)$ be the adjoint of $H(f)$.

Another notation is the dot product of two matrices. Let $A=(a_{ij})_{n\times n}$ and $B=(b_{ij})_{n\times n}$ be two $n\times n$ matrices. 
Then the dot product ``$\langle\, ,\rangle$" is defined by 
$$\langle\, A , B \, \rangle := \sum_{1\leq i,j\leq n} a_{ij}b_{ji}.$$

With these notations, we describe the Dixmier invariants of plane quartics.

Let $f,g\in K[x,y,z]_2$ be two quadratic homogeneous polynomials. 
Define
\begin{align*}
&J_{1,1}(f,g)= \langle\, H(f), H(g)\,\rangle , \\
&J_{2,2}(f,g)= \langle\, H^*(f), H^*(g)\,\rangle , \\
&J_{3,0}(f,g) = J_{3,0}(f)= \det (H(f)) ,\\
&J_{0,3}(f,g)=J_{0,3}(g) = \det(H(g)).
\end{align*}

Let $F\in K[x,y]_r$, $G\in K[x,y]_s$ be two homogeneous polynomials of degree $r$ and $s$, respectively. For $k\leq \min\{r,s\}$, define 
\begin{equation}\label{FGK}
(F,G)^k:= \frac{(r-k)!(s-k)!}{r!s!} \left.\left(\frac{\partial^2}{\partial x_1\partial y_2}- \frac{\partial^2}{\partial y_1\partial x_2}\right)^kF(x_1,y_1)G(x_2,y_2)\right|_{(x_i,y_i)=(x,y),i,1,2
}
\end{equation}

Let $P=P(x,y)\in K[x,y]_4$ be a quartic binary form. Let $Q=(P,P)^4$ defined as (\ref{FGK}). Also we let
\begin{equation}\label{SPD}\begin{split}
&\Sigma(P)=\frac{1}{2} (P,P)^4,\quad  \Psi(P)= \frac{1}{6} (P,Q)^4 \\
& \Delta(P) = \Sigma(P)^3 - 27 \Psi(P)^2\end{split}
\end{equation}
Then $\Delta(P)$ is the discriminant of $P$.

Let $u,v$ be two $K$-variables. For quartic $f\in K[x,y,z]_4$, let 
$$g=g(x,y)=f(x,y,-ux-vy).$$
Then $g(x,y)$ is a homogeneous polynomial of degree 4 with respect to the variables $x$ and $y$, and the coefficients of $g$ are expressions of $u$ and $v$. 
Thus we can define $\Sigma(g)$ and $\Psi(g)$ as in (\ref{SPD}). 
Since $\Sigma$ and $\Psi$ are expressions of the coefficients, 
we have $\Sigma(g)$ and $\Psi(g)$ are expressions of $u$ and $v$. 
An explicit computation shows that $\Sigma(g)$ and $\Psi(g)$ are polynomials of degree $2$ and $3$ in the polynomial ring $K[u,v]$ respectively.
Let $\sigma(u,v,w)$ and $\psi(u,v,w)$ be the homogenization of $\Sigma(g)$ for $w$, and $\psi(u,v,w)$ be the homogenization of $\Psi(g)$ for $w$. 
Then $\sigma(u,v,w)\in K[u,v,w]_2$ and $\psi(u,v,w)\in K[u,v,w]_3$.
Finally, we substitute $u=x, v=y, w=z$ into $\sigma(u,v,w)$ and $\psi(u,v,w)$. 
For $f\in K[x,y,z]_4$, we define
\begin{equation}\label{sp}
\begin{split}
\sigma(f)=\sigma=\sigma(x,y,z)\in K[x,y,z]_2 \\
\psi(f)=\psi=\psi(x,y,z)\in K[x,y,z]_3
\end{split}
\end{equation}
\begin{definition} Let $f\in K[x,y,z]_4$, let $\sigma$, $\psi$ defined as in (\ref{sp}). Let $\rho=D_f(\psi)$ and $\tau=D_\rho(f)$, and let ${\rm H}=\det(H(f))$.
The Dixmier invariants are defined as 
\begin{equation}
\begin{split}
&I_3=D_\sigma(f), \quad I_9= J_{1,1}(\tau, \rho), \quad I_{15}=J_{3,0}(\tau), \\
& I_6 = D_\psi ({\rm H})-8I_3^2, \quad I_{12}=J_{0,3}(\rho), \quad I_{18}= J_{2,2}(\tau, \rho) \\
& I_{27} = \Delta = \sigma^3-27\psi^2
\end{split}
\end{equation}
\end{definition}

\subsection{The Dixmier Invariants of $X_4,X_{16},X_{24}$ and $X_{96}$}\label{C369Dix}

First, we use Maxima to compute the Dixmier invariants of $X_4$. Since the equation of $X_4$ is symmetric to the parameters $r,s,u$, so are the invariants. Thus it is better to write the invariants as elementary symmetric polynomials of $r,s,u$. In the polynomial ring $K[r,s,u]$, any elementary homogeneous symmetric polynomial of $r,s,u$ of degree $d$ is uniquely determined by its leading term $r^{i_1}s^{i_2}u^{i^3}$ where $i_1\geq i_2 \geq i_3$, $i_1+i_2+i_3=d$ is an integer partition of $d$. We use $S_{[i_1,i_2,i_3]}$ to denote the symmetric polynomial whose leading term is $r^{i_1}s^{i_2}u^{i^3}$. For example, 
$$S_{[3,1,1]}=r^3su+rs^3u+rsu^3, \mbox{ and }S_{[2,1]}=r^2s+rs^2+s^2u+su^2+r^2u+ru^2.$$

\begin{proposition}\label{propX4} The Dixmier invariants of 
$$X_4(r,s,u): x^4+y^4+z^4+rx^2y^2+sy^2z^2+uz^2x^2=0$$
are

\begin{flushleft}
$
I_3=2 \left( 3 {\, S_{[2]}}+{\, S_{[1,1,1]}}+36 \right) $

\

$
I_6=62199 {\, S_{[4]}}-622086 {\, S_{[3,1,1]}}+228095 {\, S_{[2,2,2]}}+124398 {\, S_{[2,2]}}+1492776 {\, S_{[2]}}+24385464 {\, S_{[1,1,1]}}+8956656 $

\

$I_9=
81 {\, S_{[6]}}-33 {\, S_{[5,1,1]}}+15 {\, S_{[4,2,2]}}+99 {\, S_{[4,2]}}-5832 {\, S_{[4]}}+{\, S_{[3,3,3]}}+270 {\, S_{[3,3,1]}}+9936 {\, S_{[3,1,1]}}+4590 {\, S_{[2,2,2]}}+40176 {\, S_{[2,2]}}+104976 {\, S_{[2]}}+244944 {\, S_{[1,1,1]}}
$

\

$I_{12}=
-\frac{64}{729} (2025 {\, S_{[7,1,1]}}+45 {\, S_{[6,2,2]}}+9720 {\, S_{[6,2]}}-21 {\, S_{[5,3,3]}}-2997 {\, S_{[5,3,1]}}+61236 {\, S_{[5,1,1]}}-{\, S_{[4,4,4]}}-558 {\, S_{[4,4,2]}}-22032 {\, S_{[4,4]}}-47952 {\, S_{[4,2,2]}}-139968 {\, S_{[4,2]}}-7398 {\, S_{[3,3,3]}}-359640 {\, S_{[3,3,1]}}-3884112 {\, S_{[3,1,1]}}-2787696 {\, S_{[2,2,2]}}-7558272 {\, S_{[2,2]}}-34012224 {\, S_{[1,1,1]}})
$

\

$I_{15}=
-\frac{4096}{729} (243 {\, S_{[9,1,1]}}+6561 {\, S_{[8,2]}}-90 {\, S_{[7,3,3]}}-1215 {\, S_{[7,3,1]}}+255879 {\, S_{[7,1,1]}}-24 {\, S_{[6,4,4]}}-2781 {\, S_{[6,4,2]}}-15309 {\, S_{[6,4]}}-105462 {\, S_{[6,2,2]}}+1023516 {\, S_{[6,2]}}-{\, S_{[5,5,5]}}-621 {\, S_{[5,5,3]}}-14580 {\, S_{[5,5,1]}}-106110 {\, S_{[5,3,3]}}-1454355 {\, S_{[5,3,1]}}+2598156 {\, S_{[5,1,1]}}-18198 {\, S_{[4,4,4]}}-817938 {\, S_{[4,4,2]}}-4251528 {\, S_{[4,4]}}-21983724 {\, S_{[4,2,2]}}-31177872 {\, S_{[4,2]}}-8341218 {\, S_{[3,3,3]}}-85817880 {\, S_{[3,3,1]}}-410981040 {\, S_{[3,1,1]}}-463574016 {\, S_{[2,2,2]}}-510183360 {\, S_{[2,2]}}-918330048 {\, S_{[1,1,1]}})$

\

$I_{18}=
\frac{4096}{2187} (7290 {\, S_{[11,1,1]}}-2079 {\, S_{[10,2,2]}}+115911 {\, S_{[10,2]}}-1296 {\, S_{[9,3,3]}}-63342 {\, S_{[9,3,1]}}+1889568 {\, S_{[9,1,1]}}+150 {\, S_{[8,4,4]}}-12150 {\, S_{[8,4,2]}}-527796 {\, S_{[8,4]}}-1471365 {\, S_{[8,2,2]}}-1889568 {\, S_{[8,2]}}+30 {\, S_{[7,5,5]}}+7704 {\, S_{[7,5,3]}}+128628 {\, S_{[7,5,1]}}+233604 {\, S_{[7,3,3]}}-8386416 {\, S_{[7,3,1]}}-124711488 {\, S_{[7,1,1]}}+{\, S_{[6,6,6]}}+999 {\, S_{[6,6,4]}}+90882 {\, S_{[6,6,2]}}+952074 {\, S_{[6,6]}}+286686 {\, S_{[6,4,4]}}+7813422 {\, S_{[6,4,2]}}+11967264 {\, S_{[6,4]}}+32087664 {\, S_{[6,2,2]}}-289103904 {\, S_{[6,2]}}+32508 {\, S_{[5,5,5]}}+3163860 {\, S_{[5,5,3]}}+38327904 {\, S_{[5,5,1]}}+140796144 {\, S_{[5,3,3]}}+799077312 {\, S_{[5,3,1]}}+272097792 {\, S_{[5,1,1]}}+38750724 {\, S_{[4,4,4]}}+742250304 {\, S_{[4,4,2]}}+1522991808 {\, S_{[4,4]}}+8600683680 {\, S_{[4,2,2]}}+6530347008 {\, S_{[4,2]}}+4828476096 {\, S_{[3,3,3]}}+22810864896 {\, S_{[3,3,1]}}+51426482688 {\, S_{[3,1,1]}}+85710804480 {\, S_{[2,2,2]}}+33059881728 {\, S_{[2,2]}})
$

\end{flushleft}
\end{proposition}

\begin{proposition}
The Dixmier invariants of 
$$X_{16}(r,s):x^4+y^4+z^4+rx^2y^2+s(y^2z^2+z^2x^2)=0$$
are
\begin{flushleft}
$I_3=2 \left( r\, {{s}^{2}}+6 {{s}^{2}}+3 {{r}^{2}}+36\right) $,
 
 \ 
 
$I_6=\frac{1}{648}(228095 {{r}^{2}}\, {{s}^{4}}-1244172 r\, {{s}^{4}}+248796 {{s}^{4}}-622086 {{r}^{3}}\, {{s}^{2}}+248796 {{r}^{2}}\, {{s}^{2}}+24385464 r\, {{s}^{2}}+2985552 {{s}^{2}}+62199 {{r}^{4}}+1492776 {{r}^{2}}+8956656)$,

\

$I_9=\frac{64}{27} ({{r}^{3}}\, {{s}^{6}}+30 {{r}^{2}} \, {{s}^{6}}+204 r\, {{s}^{6}}+360 {{s}^{6}}+15 {{r}^{4}}\, {{s}^{4}}+540 {{r}^{3}}\, {{s}^{4}}+4788 {{r}^{2}}\, {{s}^{4}}+19872 r\, {{s}^{4}}+28512 {{s}^{4}}-33 {{r}^{5}}\, {{s}^{2}}+198 {{r}^{4}}\, {{s}^{2}}+9936 {{r}^{3}}\, {{s}^{2}}+80352 {{r}^{2}}\, {{s}^{2}}+244944 r\, {{s}^{2}}+209952 {{s}^{2}}+81 {{r}^{6}}-5832 {{r}^{4}}+104976 {{r}^{2}}),$

\

$I_{12}=\displaystyle{\frac{64 {{s}^{2}}\, {{\left( r\, {{s}^{2}}+6 {{s}^{2}}+15 {{r}^{2}}+72 r+324\right) }^{2}}\, \left( {{r}^{2}}\, {{s}^{2}}+30 r\, {{s}^{2}}+72 {{s}^{2}}-9 {{r}^{3}}+324 r\right) }{729}},$

\

$I_{15}=\displaystyle{\frac{4096}{729}} {{\left( r+3\right) }^{2}}\, {{\left( r+18\right) }^{2}}\, {{s}^{2}}\, {{\left( {{s}^{2}}+3 r+18\right) }^{2}}\, ( r\, {{s}^{4}}+6 {{s}^{4}}+18 {{r}^{2}}\, {{s}^{2}}+162 r\, {{s}^{2}}+540 {{s}^{2}}-27 {{r}^{3}}+972 r) ,$

\

$I_{18}=\frac{4096}{2187}\left( r+3\right) \, \left( r+18\right) \, {{s}^{2}}\, \left( {{s}^{2}}+3 r+18\right) \, \left( r\, {{s}^{2}}+6 {{s}^{2}}+15 {{r}^{2}}+72 r+324\right)  ({{r}^{3}}\, {{s}^{6}}+33 {{r}^{2}}\, {{s}^{6}}+228 r\, {{s}^{6}}+396 {{s}^{6}}+12 {{r}^{4}}\, {{s}^{4}}+594 {{r}^{3}}\, {{s}^{4}}+5904 {{r}^{2}}\, {{s}^{4}}+24840 r\, {{s}^{4}}+33696 {{s}^{4}}-111 {{r}^{5}}\, {{s}^{2}}-1035 {{r}^{4}}\, {{s}^{2}}+4968 {{r}^{3}}\, {{s}^{2}}+81000 {{r}^{2}}\, {{s}^{2}}+244944 r\, {{s}^{2}}+104976 {{s}^{2}}+162 {{r}^{6}}-11664 {{r}^{4}}+209952 {{r}^{2}}).$
\end{flushleft}
\end{proposition}

\begin{proposition}
The Dixmier invariants of 
$$X_{24}(r):x^4+y^4+z^4+r(x^2y^2+y^2z^2+z^2x^2)=0$$ are
\begin{flushleft}
$I_3=2 \left( {{r}^{3}}+9 {{r}^{2}}+36\right),$

\

$I_6=\frac{1}{648}(228095 {{r}^{6}}-1866258 {{r}^{5}}+559791 {{r}^{4}}+24385464 {{r}^{3}}+4478328 {{r}^{2}}+8956656)$

\

$I_9=\displaystyle{\frac{64 {{r}^{2}}\, \left( r+3\right) \, {{\left( r+18\right) }^{2}}\, {{\left( {{r}^{2}}+3 r+18\right) }^{2}}}{27}},$

\

$I_{12}=\displaystyle{\frac{64 {{r}^{3}}\, {{\left( r+18\right) }^{3}}\, {{\left( {{r}^{2}}+3 r+18\right) }^{3}}}{729}},$

\

$I_{15}=\displaystyle{\frac{4096 {{r}^{3}}\, {{\left( r+3\right) }^{3}}\, {{\left( r+18\right) }^{3}}\, {{\left( {{r}^{2}}+3 r+18\right) }^{3}}}{729}},$

\

$I_{18}=\displaystyle{\frac{4096 {{r}^{4}}\, {{\left( r+3\right) }^{2}}\, {{\left( r+18\right) }^{4}}\, {{\left( {{r}^{2}}+3 r+18\right) }^{4}}}{2187}}.$
\end{flushleft}
\end{proposition}

\begin{proposition} The Dixmier invariants of 
$$X_{96}:x^4+y^4+z^4=0$$
are
\begin{align*}
&I_3=72, \ I_6=13822, \\
&I_9=I_{12}=I_{15}=I_{18}=0.
\end{align*}
\end{proposition}

\section{The Bitangents of $X_4$,$X_{16}$,$X_{24}$ and $X_{96}$}

Explicitly, let $f=f(x,y,z)\in K[x,y,z]_4$ be the equation of a plane quartic $C$. 
Let $L: ax+by+cz=0$, $a,b,c\in K$ be a line in ${\mathbb P}^2_{(x,y,z)}$. Thus the point $(a,b,c)\in {\mathbb P}^2_{(a,b,c)}$ determines the line $L$. 
So without lost of generality, we can assume that $c\neq 0$, and say $c=1$. 
This time $L: ax+by+z=0$ gives the condition $z=-ax-by$. 
Substitute this relation into $f(x,y,z)$ we have a quadratic form $f(x,y,-ax-by)\in R[x,y,z]_2$ where $R=K[a,b]$. 
If $L$ is a bitangent for some $a,b\in K$, then there exist $\lambda_0,\lambda_1,\lambda_2\in K$ such that 
\begin{equation}\label{lambda012}f(x,y,-ax-by)=(\lambda_0x^2+\lambda_1xy+\lambda_2y^2)^2.\end{equation}
\begin{definition}\label{IJ}For any quartic $f\in K[x,y,z]_4$, let $I(f)$ be the ideal of $K[a,b,\lambda_0,\lambda_1,\lambda_2]$ generated by comparing the coefficients of both sides of the monomials of  $x,y$ in the expansion of (\ref{lambda012}). 
Let $J(f)$ be elimination ideal of $I$ with respect to $\lambda_0,\lambda_1,\lambda_2$ in $K[a,b]$. \end{definition}  
The ideal $J(f)$ gives the conditions of $L$ being a bitangent of $C$. In general one cannot solve $a,b$ over ${\mathbb Q}$, and even there exists $L$ such that $a,b\in {\mathbb Q}$, the tangency points $p_1,p_2$ are not $\mathbb Q$-rational points of $C$.

\begin{theorem}\label{BitanX4}
The curve $X_4$ has symbolic solutions for all the 28 bitangents.
\end{theorem}

{\it Proof}\qquad Let $f$ be the equation defined by
\begin{equation}\label{X4}
f(x,y,z)=x^4+y^4+z^4+rx^2y^2+sy^2z^2+uz^2x^2=0.
\end{equation}
Let $J(f)$ be the ideal defined by Definition \ref{IJ}. Using Macaulay2, we can compute the primary decomposition of $J(f)$. We can input
\begin{verbatim}
R = QQ[r,s,u,a,b,k_0,k_1,k_2][x,y,z]
f = x^4+y^4+z^4+r*x^2*y^2+s*y^2*z^2+u*z^2*x^2
g = (k_0*x^2+k_1*x*y+k_2*y^2)^2
h = substitute(f,{z => -a*x-b*y})
H= h-g
Coe = coefficients H
L = flatten entries Coe#1
S = QQ[r,s,u,a,b,k_0,k_1,k_2]
I = ideal L
psi=map(S,R)
phi=map(R,S)
J = psi I
E=eliminate(J,{k_0,k_1,k_2})
T = QQ[r,s,u,a,b]
xi=map(T,S)
U = xi E
D = primaryDecomposition U
\end{verbatim}
The output says $J(f)$ has 3 irreducible components. They are

\

$J_1=\langle s^{2}a^{4}-u^{2}b^{4}+s^{2}u\,a^{2}-s\,u^{2}b^{2}-4\,a^{4}+4\,b^{4}-4\,u\,a^{2}+4\,s\,b^{2}+s^{2}-u^{2},\quad \\ r\,u^{2}a^{2}b^{2}-r^{2}a^{4}+u^{2}a^{4}+u^{2}b^{4}-r^{2}u\,a^{2}-4\,r\,a^{2}b^{2}-4\,b^{4}+4\,u\,a^{2}-r^{2}+4,\quad \\ r\,s^{2}a^{2}b^{2}-r^{2}b^{4}+s^{2}b^{4}+u^{2}b^{4}-s^{2}u\,a^{2}-r^{2}s\,b^{2}+s\,u^{2}b^{2}-4\,r\,a^{2}b^{2}-4\,b^{4}+4\,u\,a^{2}-r^{2}-s^{2}+u^{2}+4,\quad \\ s\,u\,a^{2}b^{4}+u^{2}b^{6}+s\,u^{2}b^{4}-2\,r\,a^{2}b^{4}-2\,r\,s\,a^{2}b^{2}-r\,u\,b^{4}-4\,b^{6}-r\,s\,u\,b^{2}+4\,u\,a^{2}b^{2}-2\,s\,b^{4}+s\,u\,a^{2}+u^{2}b^{2}-2\,r\,a^{2}-r\,u+4\,b^{2}+2\,s,\quad \\ s^{2}a^{2}b^{4}+s\,u\,b^{6}+s^{2}u\,b^{4}-2\,r\,b^{6}-3\,r\,s\,b^{4}-4\,a^{2}b^{4}-r\,s^{2}b^{2}+2\,u\,b^{4}-s^{2}a^{2}+3\,s\,u\,b^{2}-2\,r\,b^{2}-r\,s+4\,a^{2}+2\,u,\quad \\r\,s\,a^{2}b^{4}-r\,u\,b^{6}-r\,s\,u\,b^{4}-2\,u\,a^{2}b^{4}+2\,s\,b^{6}-2\,s\,u\,a^{2}b^{2}+r^{2}b^{4}+r^{2}s\,b^{2}+4\,r\,a^{2}b^{2}+r\,s\,a^{2}-r\,u\,b^{2}+4\,b^{4}-2\,u\,a^{2}-2\,s\,b^{2}+r^{2}-4,\quad \\ s\,u\,a^{4}b^{2}+u^{2}a^{2}b^{4}+s\,u^{2}a^{2}b^{2}-2\,r\,a^{4}b^{2}-r\,s\,a^{4}-2\,r\,u\,a^{2}b^{2}-4\,a^{2}b^{4}-r\,s\,u\,a^{2}+2\,u\,a^{4}+u^{2}a^{2}+s\,u\,b^{2}-2\,r\,b^{2}-r\,s+4\,a^{2}+2\,u,\quad \\ r\,s\,a^{4}b^{2}-r\,u\,a^{2}b^{4}-2\,u\,a^{4}b^{2}+2\,s\,a^{2}b^{4}-s\,u\,a^{4}+s\,u\,b^{4}+2\,r\,a^{4}-2\,r\,b^{4}+r\,u\,a^{2}-r\,s\,b^{2}-2\,s\,a^{2}+2\,u\,b^{2},\quad \\ 
s\,u\,a^{6}+u^{2}a^{4}b^{2}+s\,u^{2}a^{4}-2\,r\,a^{6}-3\,r\,u\,a^{4}-4\,a^{4}b^{2}-r\,u^{2}a^{2}+2\,s\,a^{4}+3\,s\,u\,a^{2}-u^{2}b^{2}-2\,r\,a^{2}-r\,u+4\,b^{2}+2\,s,\quad \\ r\,s\,a^{6}-r\,u\,a^{4}b^{2}+r\,s\,u\,a^{4}-2\,u\,a^{6}+2\,s\,a^{4}b^{2}-r^{2}a^{4}+2\,s\,u\,a^{2}b^{2}-r^{2}u\,a^{2}-4\,r\,a^{2}b^{2}+r\,s\,a^{2}-4\,a^{4}-r\,u\,b^{2}+2\,u\,a^{2}+2\,s\,b^{2}-r^{2}+4\rangle$
   
      \
     
      $J_2=\langle a,u^{2}b^{4}+2\,r\,u\,b^{2}-4\,b^{4}-4\,s\,b^{2}+r^{2}-4\rangle$
      
      \
      
      $J_3=\langle b,s^{2}a^{4}+2\,r\,s\,a^{2}-4\,a^{4}-4\,u\,a^{2}+r^{2}-4\rangle$

\

For $J_2$, for example, we have $a=0$, and $b$ satisfies an equation of degree 4, which is solvable. Similarly to $J_3$, we have  another 4 bitangents. Thus we have 8 bitangents. Remember that $f$ is symmetric with respect to $r,s$ and $u$, and we are considering the affine $xy$-plane of the projective plane ${\mathbb P}^2$. Thus if we consider $yz$-plane and the $zx$-plane, we have another such ideals $J_2'$, $J_2''$, $J_3'$, $J_3''$, each ideal gives 4 solutions of bitangents. 

For $J_2$, the condition $a=0$ implies that the equation of the bitangent looks like $by+z=0$. This line does not lie on the $yz$-plane. On the other hand, if we consider the $yz$-plane, let the line has equation $x+by+cz=0$, and run the same algorithm, then we have another ideal
$$J_2''=\langle c,u^{2}b^{4}+2\,s\,u\,b^{2}-4\,b^{4}-4\,r\,b^{2}+s^{2}-4\rangle$$
the condition $c=0$ implies that the 4 bitangents given by $J_2''$ do not lie on the $xy$-plane. The algorithm will output all the bitangents on the corresponding affine plane, so $J_1,J_2$ and $J_3$ will solve 24 bitangents because there are 28 inall. Thus the ideal $J_1$ gives $28-8-4=16$ bitangents.

Last we consider $J_1$. Eiminate $J_1$ with respect to $a$, we have that $b$ satsifies the degree 8 equation
\begin{flushleft}
$-{{b}^{4}}\, {{s}^{2}}\, {{u}^{2}}-2 {{b}^{6}} s\, {{u}^{2}}-2 {{b}^{2}} s\, {{u}^{2}}-{{b}^{8}}\, {{u}^{2}}-2 {{b}^{4}}\, {{u}^{2}}-{{u}^{2}}+{{b}^{6}} r\, {{s}^{2}} u+{{b}^{2}} r\, {{s}^{2}} u+{{b}^{8}} r s u+6 {{b}^{4}} r s u+r s u+4 {{b}^{6}} r u+4 {{b}^{2}} r u-{{b}^{4}}\, {{r}^{2}}\, {{s}^{2}}-{{b}^{8}}\, {{s}^{2}}+2 {{b}^{4}}\, {{s}^{2}}-{{s}^{2}}-2 {{b}^{6}}\, {{r}^{2}} s-2 {{b}^{2}}\, {{r}^{2}} s-{{b}^{8}}\, {{r}^{2}}-2 {{b}^{4}}\, {{r}^{2}}-{{r}^{2}}+4 {{b}^{8}}-8 {{b}^{4}}+4=0.$
\end{flushleft}
This equation contains only even power terms of $b$, so let $B=b^2$, then we have a degree 4 equation

\begin{flushleft}
${{B}^{2}}\, \left( -{{s}^{2}}\, {{u}^{2}}-2 {{u}^{2}}+6 r s u-{{r}^{2}}\, {{s}^{2}}+2 {{s}^{2}}-2 {{r}^{2}}-8\right) +{{B}^{3}}\, \left( -2 s\, {{u}^{2}}+r\, {{s}^{2}} u+4 r u-2 {{r}^{2}} s\right) +B\, \left( -2 s\, {{u}^{2}}+r\, {{s}^{2}} u+4 r u-2 {{r}^{2}} s\right) -{{u}^{2}}+{{B}^{4}}\, \left( -{{u}^{2}}+r s u-{{s}^{2}}-{{r}^{2}}+4\right) +r s u-{{s}^{2}}-{{r}^{2}}+4=0$
\end{flushleft}
of $B$, which is solvable. Observe that the lowest degree of the generators in $J_1$ for $a$ is 2 (for example the third generator), so for each fixed $b$ one can solve a pair of $a$'s, and 8 $b$'s give 16 bitangents.\qquad $\blacksquare$

For$$X_{16}(r,s):f(x,y,z)=x^4+y^4+z^4+rx^2y^2+s(y^2z^2+z^2x^2)=0,$$
the primary decomposition of $J(f)$ is
$J_1=\langle a,s^{2}b^{4}+2\,r\,s\,b^{2}-4\,b^{4}-4\,s\,b^{2}+r^{2}-4\rangle,\,$ 
      
      $J_2=\langle b,s^{2}a^{4}+2\,r\,s\,a^{2}-4\,a^{4}-4\,s\,a^{2}+r^{2}-4\rangle,\,$
      
      $J_3=\langle s+2,r-2\rangle,\,$
      
      $J_4=\langle s-2,r-2\rangle,\,$
      
      $J_5=\langle a+b,s^{2}b^{4}-r\,b^{4}-r\,s\,b^{2}-2\,b^{4}+2\,s\,b^{2}-r+2\rangle,\,$
      
      $J_6=\langle a-b,s^{2}b^{4}-r\,b^{4}-r\,s\,b^{2}-2\,b^{4}+2\,s\,b^{2}-r+2\rangle,\,$
      
      $J_7=\langle a^{2}+b^{2}+s,r\,b^{4}+r\,s\,b^{2}-2\,b^{4}-2\,s\,b^{2}-s^{2}+r+2\rangle$

The component $J_4$ gives $r=s=2$, which is the case $X=X_{24}$. For $J_3$, we have $r=-s=2$. Let $\zeta_4$ be the primitive 4th root of $1$, then $x\mapsto \zeta_4x, y\mapsto \zeta_4 y, z\mapsto z$ is a projective isomorphism making $r\mapsto -r$, so this also gives the situation when $X=X_{24}$. Otherwise, the second generator for each of the components are quartic polynomials with one variable such that only even degree terms occure. Thus they are essentially quadratic equations. For example, the first component $J_1$ gives four bitangents, they are 
$by+z$ where
\begin{align*}&b=-\sqrt{\frac{2 \sqrt{\left( 2-r\right) \, {{s}^{2}}+{{r}^{2}}-4}+\left( 2-r\right)  s}{{{s}^{2}}-4}}, \\ &b=\sqrt{\frac{2 \sqrt{\left( 2-r\right) \, {{s}^{2}}+{{r}^{2}}-4}+\left( 2-r\right)  s}{{{s}^{2}}-4}}, \\ &b=-\sqrt{-\frac{2 \sqrt{\left( 2-r\right) \, {{s}^{2}}+{{r}^{2}}-4}+\left( r-2\right)  s}{{{s}^{2}}-4}}, \\ &b=\sqrt{-\frac{2 \sqrt{\left( 2-r\right) \, {{s}^{2}}+{{r}^{2}}-4}+\left( r-2\right)  s}{{{s}^{2}}-4}}.
\end{align*}

For 
$$X_{24}(r):f(x,y,z)=x^4+y^4+z^4+r(x^2y^2+y^2z^2+z^2x^2)=0,$$
the primary decomposition of $J(f)$ is

$
J_1=\langle r-2\rangle,\,$

$J_2=\langle a,r\,b^{4}+2\,b^{4}+2\,r\,b^{2}+r+2\rangle,\,$

$J_3=\langle b,r\,a^{4}+2\,a^{4}+2\,r\,a^{2}+r+2\rangle,\,$

$J_4=\langle b+1,a-1\rangle,\,$

$J_5=\langle b-1,a+1\rangle,\,$

$J_6=\langle a+b,r\,b^{2}+b^{2}+1\rangle,\,$
 
$J_7=\langle b+1,a+1\rangle,\,$

$J_8=\langle b-1,a-1\rangle,\,$

$J_9=\langle a-b,r\,b^{2}+b^{2}+1\rangle,\,$

$J_{10}=\langle a+1,b^{2}+r+1\rangle,\,$

$J_{11}=\langle a-1,b^{2}+r+1\rangle,\,$

$J_{12}=\langle b+1,a^{2}+r+1\rangle,\,$

$J_{13}=\langle b-1,a^{2}+r+1\rangle.
$

Thus, it is easy to write down the 28 bitangents of $X_{24}$ for $r\neq 2$.

\begin{remark} The case $r=2$, and furthermore, $|r|=|s|=|u|=2$ correspond to the situation when $X$  degenerates to a total square of a quadric polynomial, or say, a double conic.
\end{remark}

In the end, for $$X_{96}:x^4+y^4+z^4=0,$$
there are 16 bitangents $ax+by+cz=0$
where
$$a^2,b^2,c^2=\pm 1,$$
and 12 bitangents given by one of $a,b,c$ is 0, another is 1, and the one left is a 4th root of $-1$.

\section{The Matrix Representation Problem}

Let $f(x,y,z)\in \kappa [x,y,z]_4$ be a homogeneous polynomial of degree $4$ over the algebraic closed field $\kappa$, and let $$X:f(x,y,z)=0$$ be the plane quartic curve defined by $f$.

The matrix representation problem for $X$ asks whether the polynomial $f(x,y,z)$  could be written of the form
$$f(x,y,z)={\rm det}(xA+yB+zC)$$
where $A,B,C$ are symmetric matrices whose entried defined over $\kappa$. According to Section 2 in  \cite{Sturmfels2}, if 
\begin{equation}\label{condition}f(x,0,0)=x^4\quad \mbox{and}\quad f(x,y,0)=\prod_{i=1}^4 (x+\beta_iy) 
\end{equation}
for some $\beta_1,\beta_2,\beta_3,\beta_4\in \kappa$, then one can assume that
$$
 A =\begin{pmatrix}1 &&& \\
 &1 && \\ && 1 & \\ &&& 1 \end{pmatrix}, \quad B = \begin{pmatrix}\beta_1 &&& \\
 & \beta_2 && \\ && \beta_3 & \\ &&& \beta_4 \end{pmatrix} 
,\quad C=\begin{pmatrix}{c_{11}} & {c_{12}} & {c_{13}} & {c_{14}}\\
{c_{12}} & {c_{22}} & {c_{23}} & {c_{24}}\\
{c_{13}} & {c_{23}} & {c_{33}} & {c_{34}}\\
{c_{14}} & {c_{24}} & {c_{34}} & {c_{44}}\end{pmatrix}.$$

Furthermore, we have
\begin{equation}\label{cii}c_{ii}=\beta_i\cdot \frac{\frac{\partial f}{\partial z}(-\beta_i,1,0)}{\frac{\partial f}{\partial y}(-\beta_i,1,0)},\quad i=1,2,3,4.
\end{equation}

Let $X=X_4$ defined by 
\begin{equation}\label{X4}f(x,y,z)= x^4+y^4+z^4+rx^2y^2+sy^2z^2+uz^2x^2=0.
\end{equation}
We check the condtions in (\ref{condition}). Obviously $f(x,0,0)=x^4$. We also have
$$ f(x,y,0)= x^4+rx^2y^2+y^4 
= (x+py)(x-py)(x+qy)(x-qy) $$
where 
\begin{equation}\label{pq} p= \frac{\sqrt{-\sqrt{{{r}^{2}}-4}-r}}{\sqrt{2}}, q=\frac{\sqrt{\sqrt{{{r}^{2}}-4}-r}}{\sqrt{2}}\, .
\end{equation}
Hence 
$$B=\begin{pmatrix}p & 0 & 0 & 0\\
0 & -p & 0 & 0\\
0 & 0 & q & 0\\
0 & 0 & 0 & -q\end{pmatrix}.$$

Next, the partial derivative $$\frac{\partial f}{\partial z}=-3yz^2$$
 implies that if $z=0$, so $c_{ii}=0$ for $i=1,2,3,4$ by (\ref{cii}).

For convinience we denote
$$D=\begin{pmatrix} & {c_{12}} & {c_{13}} & {c_{14}}\\
& & {c_{23}} & {c_{24}}\\
&& & {c_{34}}\\
&&& \end{pmatrix} = \begin{pmatrix} & a & b & d\\
& & c & e\\
&& & f\\
&&& \end{pmatrix},$$
then $C=D+\,^{\sf t}D$ where $\,^{\sf t}D$ is the matrix transpose of $D$ since $c_{ii}=0$ for $i=1,2,3,4$.

Using Maxima, we directly compute the coefficients of $$\det (xA+yB+zC)=\det \begin{pmatrix} x+py & a z & b z & d z\\
a z & x-py & c z & e z\\
b z & c z & x+q y & f z\\
d z & e z & f z & x-q y\end{pmatrix}$$ and compare the coefficients with $f(x,y,z)$ in (\ref{X4}), the output is a system of equations
\begin{align}
& -{{e}^{2}} q-{{d}^{2}} q+{{c}^{2}} q+{{b}^{2}} q-{{e}^{2}} p+{{d}^{2}} p-{{c}^{2}} p+{{b}^{2}} p=0, \label{oeq1} \\
& -s+{{a}^{2}}\, {{q}^{2}}-{{e}^{2}} p q+{{d}^{2}} p q+{{c}^{2}} p q-{{b}^{2}} p q+{{f}^{2}}\, {{p}^{2}}=0,\label{oeq2}\\ 
&-u-{{f}^{2}}-{{e}^{2}}-{{d}^{2}}-{{c}^{2}}-{{b}^{2}}-{{a}^{2}}=0,\label{oeq3}\\
& 2 a d e q-2 a b c q+2 c e f p-2 b d f p=0, \label{oeq4}\\
& 2 c e f+2 b d f+2 a d e+2 a b c=0,\label{oeq5}\\
& {{a}^{2}}\, {{f}^{2}}-2 a b e f-2 a c d f+{{b}^{2}}\, {{e}^{2}}-2 b c d e+{{c}^{2}}\, {{d}^{2}}-1=0,\label{oeq6}\\
& -r-{{q}^{2}}-{{p}^{2}}=0, \label{oeq7}\\
& {{p}^{2}}\, {{q}^{2}}-1 \label{oeq8}.
\end{align}

The last two equations (\ref{oeq7}) and (\ref{oeq8}) are identities. According to (\ref{pq}), 
$pq=1$. 
We rewrite the equation system as
\begin{align}
& (b^2-e^2)(q+p)+(c^2-d^2)(q-p) =0, \label{e1}\\ 
& a^2q^2-b^2+c^2+d^2-e^2+f^2p^2=s \label{e2}\\
&{{f}^{2}}+{{e}^{2}}+{{d}^{2}}+{{c}^{2}}+{{b}^{2}}+{{a}^{2}}=-u,\label{e3}\\
&  a d e q- a b c q+ c e f p- b d f p=0=0,\label{e4}\\
&  c e f+ b d f+ a d e+ a b c=0,\label{e5}\\
& {{a}^{2}}\, {{f}^{2}}-2 a b e f-2 a c d f+{{b}^{2}}\, {{e}^{2}}-2 b c d e+{{c}^{2}}\, {{d}^{2}}=1\label{e6}.
\end{align}
of the 6 variables $a,b,c,d,e,f$.

\begin{theorem}\label{MRPX4} The equation system (\ref{e1})-(\ref{e6}) has a symbolic solution.
\end{theorem}

{\it Proof}\qquad 

Generically, it does not matter which parameters we wish to eliminate first. However, if we wish to eiliminate $a,f$ using (\ref{e4}) and (\ref{e5}), then these are linear equations and  $a=f=0$. Then the equation system becomes

\begin{align}
& (b^2-e^2)(q+p)+(c^2-d^2)(q-p) =0, \label{eqq1}\\ 
& -(b^2+e^2)+(c^2+d^2)=s \label{eqq2}\\
&(b^2+e^2)+(c^2+d^2)=-u,\label{eqq3}\\
& (be-cd)^2=1\label{eqq6}.
\end{align}

We only seek for one solution to the equation system , thus if there is an "either-or" argument in any step, we can choose one of them as our solution. For example, we can choose 
\begin{equation}\label{eqq7}
be=cd+1
\end{equation}
in (\ref{eqq6}). From (\ref{eqq2}) and (\ref{eqq3}) we have
\begin{align}
& b^2+e^2 = -\frac{u+s}{2} \label{eqq9} \\
& c^2+d^2 = -\frac{u-s}{2} \label{eqq10}
\end{align}

 Write (\ref{eqq1})$^2$ as
\begin{equation}\label{eqq8}
(b^2+e^2+2be)(b^2+e^2-2be)(q+p)^2=(c^2+d^2+2cd)(c^2+d^2-2cd)(q-p)^2
\end{equation}
then substitude (\ref{eqq7}),(\ref{eqq8}) and (\ref{eqq10}) into (\ref{eqq8}), we have
\begin{equation}\label{eq35}
\left(\frac{p-q}{p+q}\right)^2\left[\frac{(u-s)^2}{4}-4(cd)^2\right]+4(cd+1)^2-\frac{(u+s)^2}{4}=0.
\end{equation}
Last, from (\ref{oeq7}) and (\ref{oeq8}) we have
$$\left(\frac{p-q}{p+q}\right)^2= \frac{r+2}{r-2}$$
so (\ref{eq35}) becomes
\begin{equation}\left(\label{eq36}
\frac{r+2}{r-2}\right)\left[\frac{(u-s)^2}{4}-4(cd)^2\right]+4(cd+1)^2-\frac{(u+s)^2}{4}=0.
\end{equation}
which is a quadrtic equation of $cd$. Together with (\ref{eqq10}) we have $c$ and $d$, respectively, and thus (\ref{eqq9}) and (\ref{eqq7}) gives $b$ and $e$. \qquad $\blacksquare$

\begin{flushright}

{\sc School of Mathematics

Sun Yat-Sen University

Guangzhou China

510275
}

{\tt liangdun@mail.sysu.edu.cn}
\end{flushright}

\end{document}